\documentclass[a4paper,11pt]{article}
\usepackage{amsmath}
\usepackage{amsfonts}
\usepackage{amsthm}
\usepackage[round]{natbib}
\usepackage{url}

\begin{document}
\bibliographystyle{plainnat}

\begingroup
\newtheorem{thm}{Theorem}
\newtheorem{cor}[thm]{Corollary}      
\newtheorem{lem}[thm]{Lemma}        
\newtheorem{prop}[thm]{Proposition} 
\newtheorem{conj}[thm]{Conjecture}  
\endgroup

\newtheorem{construction}{Construction}
\newtheorem{exmp}{Example}       
\newtheorem{defn}{Definition}  
\newtheorem{rem}{Remark}  
\newtheorem{notation}{Notation.}
\renewcommand{\thenotation}{}
\newtheorem{terminology}{Terminology.}
\renewcommand{\theterminology}{}


\title{Uniform semi-Latin squares and their pairwise-variance aberrations}

\author{R. A. Bailey\\
School of Mathematics and Statistics\\
University of St Andrews\\
North Haugh, St Andrews, Fife KY16 9SS, UK\\
{\small\texttt{rab24@st-andrews.ac.uk}}\\
\and\\
Leonard H. Soicher\\
School of Mathematical Sciences\\
Queen Mary University of London\\
Mile End Road, London E1 4NS, UK\\
{\small\texttt{L.H.Soicher@qmul.ac.uk}}\\
}

\maketitle

\begin{abstract}
\parindent=0pt
\parskip=3pt

For integers $n>2$ and $k>0$, an $(n\times n)/k$ 
\textit{semi-Latin square} is 
an $n\times n$ array of $k$-subsets (called \textit{blocks}) of an 
$nk$-set (of \textit{treatments}), such that each treatment occurs 
once in each row and once in each column of the array.
A semi-Latin square is \textit{uniform} if every pair of blocks, 
not in the same row or column, intersect in the same positive number 
of treatments. 
It is known that a uniform $(n\times n)/k$ semi-Latin 
square is Schur optimal in
the class of all $(n\times n)/k$ semi-Latin squares, and here
we show that when a uniform $(n\times n)/k$ semi-Latin square exists, 
the Schur optimal $(n\times n)/k$ semi-Latin squares are
precisely the uniform ones. 
We then compare uniform semi-Latin squares using the criterion of
pairwise-variance (PV) aberration, introduced by J.~P. Morgan for affine
resolvable designs, and determine the uniform $(n\times n)/k$ 
semi-Latin squares with minimum PV aberration when there
exist $n-1$ mutually orthogonal Latin squares of order $n$. 
These do not exist when $n=6$, and the smallest 
uniform semi-Latin squares in this case have size $(6\times 6)/10$.
We present a complete classification of the uniform $(6\times 6)/10$ 
semi-Latin squares, 
and display the one with least PV aberration. 
We give a construction producing a
uniform $((n+1)\times (n+1))/((n-2)n)$ semi-Latin square 
when there exist $n-1$ mutually orthogonal Latin squares of order $n$,
and determine the PV aberration of such a uniform semi-Latin square.
Finally, we describe how certain affine resolvable designs and
balanced incomplete-block designs can be constructed from uniform
semi-Latin squares. From the
uniform $(6\times 6)/10$ semi-Latin squares we classified, we obtain 
(up to block design isomorphism) exactly
16875 affine resolvable 
designs for 72 treatments in 36 blocks of
size~12 and 8615 balanced incomplete-block designs 
for 36 treatments in 84 blocks of size~6.   
In particular, this shows that there are at least 16875 
pairwise non-isomorphic orthogonal arrays $\mathrm{OA}(72,6,6,2)$.
\medskip

\medskip

\end{abstract}

\section{Introduction}
For integers $n>2$ and $k>0$, an $(n\times n)/k$ \textit{semi-Latin square} is 
an $n\times n$ array of $k$-subsets (called \textit{blocks}) of an 
$nk$-set (of \textit{treatments}), such that each treatment occurs once in
each row and once in each column of the array. Note that an $(n\times n)/1$ semi-Latin square is 
the same thing as a Latin square of order $n$. 
We consider two $(n\times n)/k$ semi-Latin squares to be \textit{isomorphic} 
if one can be obtained from the other by applying an \textit{isomorphism}, 
which is a sequence of zero or more of: permuting the rows, permuting the columns, 
transposing the array, and renaming the treatments. An \textit{automorphism} of a semi-Latin square 
$S$ is an isomorphism mapping $S$ onto itself.  
The applications of semi-Latin squares include the design of agricultural experiments, 
consumer testing, and via their duals, human-machine interaction (see \citet{Bailey92,Bailey11}).

A $(v,b,r,k)$-\textit{design} is a binary block design 
for $v$ treatments in $b$ blocks of size $k$ (considered as $k$-subsets 
of the set of treatments), such that each treatment 
is in exactly $r$ blocks.
If we ignore the block structure of an $(n\times n)/k$ semi-Latin 
square $S$ then we obtain an $(nk,n^2,n,k)$-design called 
the \textit{underlying block design} of $S$. 
A $(v,b,r,k)$-design with $k<v$ and $r>1$ is \textit{resolvable} if its collection of 
blocks can be partitioned
into $r$ partitions of the treatments (called \textit{parallel classes}),
and such a resolvable design is \textit{affine resolvable} if every pair of blocks
in distinct parallel classes intersect in the same positive number $\mu$ of treatments.
A $(v,b,r,k)$-design is a $(v,b,r,k,\lambda)$-\textit{balanced
incomplete-block design} (\textit{BIBD}) if $1<k<v$ and every pair
of distinct treatments occur together in exactly $\lambda$ blocks.
Two $(v,b,r,k)$-designs are \textit{isomorphic} (as block designs)
if there is a bijection from the treatments of the first to those of
the second which maps the list of blocks of the first onto that of
the second in some order. Such a bijection is called a (block design)
\textit{isomorphism}, and an \textit{automorphism} of a $(v,b,r,k)$-design
is an isomorphism from that block design to itself.

An \textit{orthogonal array} of \textit{strength} $t$ with $N$ rows, $r$
columns ($r\ge t$), and based on $s$ symbols ($s\ge 2$), here taken to be
$1,2,\ldots s$, or an \textit{orthogonal array} $\mathrm{OA}(N,r,s,t)$,
is an $N\times r$ array whose entries are symbols, such that
for every $N\times t$ subarray, each of the possible
$s^t$ $t$-tuples of symbols occurs as a row equally often 
(which must be $N/s^t$ times).
As \cite{HSS99} point out, there are many trivial constructions of
orthogonal arrays of strength one, so we ignore this case.
Two orthogonal arrays $\mathrm{OA}(N,r,s,t)$ are \textit{isomorphic}
if one can be obtained from the other by permuting the rows, permuting
the columns, and permuting the symbols separately within each column. It
is known that orthogonal arrays of strength~$2$ and affine resolvable
designs are equivalent combinatorial objects. In particular, \citet{BMM95}
describe how to construct an equivalent $\mathrm{OA}(v,r,s,2)$ from an
affine resolvable $(v,b,r,v/s)$-design (see also \citet{Morgan10}),
such that two affine resolvable $(v,b,r,v/s)$-designs are isomorphic
if and only if their equivalent orthogonal arrays are isomorphic. We
remark that affine resolvable designs were introduced by \citet{Bose42}
(in the context of BIBDs), and orthogonal arrays were introduced by
\citet{RRao47}.

An $(n\times n)/k$ semi-Latin square $S$ is \textit{uniform} if every pair of blocks, 
not in the same row or column, intersect in the same positive number 
$\mu=\mu(S)$ of treatments (in which case $k=\mu(n-1)$).
For example, here is a $(3\times 3)/4$ uniform
semi-Latin square with $\mu=2$: 
\begin{equation}
\label{USLSexample}
\begin{tabular}{|r|r|r|}
\hline
1 4 7 10 &
2 5 8 11 &
3 6 9 12 \\
\hline
3 6 8 11 &
1 4 9 12 &
2 5 7 10 \\
\hline
2 5 9 12 &
3 6 7 10 &
1 4 8 11 \\
\hline
\end{tabular}
\end{equation}

Uniform semi-Latin squares were introduced, constructed, and studied
by \citet{Soicher12}, where it was shown
that a uniform $(n\times n)/k$ semi-Latin square is Schur optimal 
(defined later) in
the class of all $(n\times n)/k$ semi-Latin squares.  

In this paper, we further the study of uniform semi-Latin squares.  
We show that, if a uniform $(n\times n)/k$ semi-Latin square exists, 
then the Schur optimal $(n\times n)/k$ semi-Latin squares are
precisely the uniform ones. 

We then compare uniform $(n\times n)/k$ semi-Latin 
squares using the criterion 
of pairwise-variance (PV) aberration, introduced by \citet{Morgan10} for
affine resolvable designs, and determine 
the uniform $(n\times n)/k$ semi-Latin squares with minimum
PV aberration when there
exist $n-1$ mutually orthogonal Latin squares (MOLS) of order $n$.
These do not exist when $n=6$, and the smallest 
uniform semi-Latin squares in this case have size $(6\times 6)/10$.
We describe a complete classification of the
uniform $(6\times 6)/10$ semi-Latin squares, 
and find that, up to isomorphism, there
are exactly 8615 such designs. We compare their PV
aberrations, and display the one with least PV aberration. 

We give a construction producing a
uniform $((n+1)\times (n+1))/((n-2)n)$ semi-Latin square 
when there exist $n-1$ MOLS of order $n$,
and determine the PV aberration of such a uniform semi-Latin square.

Finally, we describe how a uniform $(n\times n)/(\mu(n-1))$ semi-Latin
square can be used to construct 
two (possibly isomorphic) affine resolvable $(\mu n^2,n^2,n,\mu n)$-designs 
and an $(n^2, \mu n(n+1),\mu (n+1),n,\mu)$-BIBD. From the
uniform $(6\times 6)/10$ semi-Latin squares we classified, we obtain 
(up to block design isomorphism) exactly
16875 affine resolvable $(72,36,6,12)$-designs
and 8615 $(36,84,14,6,2)$-BIBDs.
In particular, this shows that there are at least 
16875 pairwise non-isomorphic orthogonal arrarys $\mathrm{OA}(72,6,6,2)$.

\section{Schur optimality}

Let $\Delta$ be a $(v,b,r,k)$-design.
The \textit{concurrence matrix}
$\Lambda$ of $\Delta$ is the $v\times v$ matrix whose rows
and columns are indexed by the treatments of $\Delta$, and whose
$(\alpha,\beta)$-entry is the number of blocks containing both
$\alpha$ and $\beta$  (this entry is the
\textit{concurrence} of treatments $\alpha$ and $\beta$).

The \textit{scaled information matrix} of $\Delta$ is 
\[ F(\Delta):=I_v-(rk)^{-1}\Lambda,\]
where $I_v$ denotes the $v\times v$ identity matrix. 
The eigenvalues of $F(\Delta)$
are all real and lie in the interval $[0,1]$.
The all-$1$ vector is an eigenvector of $F(\Delta)$ with corresponding
eigenvalue $0$. The remaining eigenvalues (counting repeats) 
are called the \textit{canonical efficiency factors} of $\Delta$. It is well
known that these canonical efficiency factors are all non-zero if and
only if $\Delta$ is connected (that is, its treatment-block incidence
graph is connected), and they are all equal to $1$ if and only if $k=v$.

Now suppose that $\Delta$ has canonical efficiency factors
$\delta_1\le \cdots\le \delta_{v-1}$.  We say that $\Delta$ is \textit{Schur
optimal} in a class $\cal C$ of $(v,b,r,k)$-designs containing
$\Delta$ if for each design $\Gamma\in {\cal C}$, with canonical
efficiency factors $\gamma_1\le \cdots \le\gamma_{v-1}$, we have
\[ \sum_{i=1}^{\ell}\delta_i\ge \sum_{i=1}^{\ell}\gamma_i,\] for
$\ell=1,\ldots,v-1$.  A Schur optimal design need not exist within a
given class $\cal C$, but, when it does, that design is optimal in $\cal C$ 
with respect to a very wide range of statistical optimality criteria,
including being $\Phi_p$-optimal, for all $p\in (0,\infty)$, and also A-,
D-, and E-optimal (see \citet{GW81}; see also \citet{BC09} or \citet{SS89}
for definitions of these optimality criteria).

As recommmended by \citet{Bailey92}, for the purposes of statistical optimality,
we compare $(n\times n)/k$ semi-Latin
squares as their underlying $(nk,n^2,n,k)$-designs. Thus, we take 
the \textit{canonical efficiency factors} of a semi-Latin square to be those 
of its underlying block design, and to say that an $(n\times n)/k$ semi-Latin
square $S$ is \textit{Schur optimal} means that its underlying block design
is Schur optimal in the class of underlying block designs of
$(n\times n)/k$ semi-Latin squares. 

The \textit{dual} of the $(v,b,r,k)$-design
$\Delta$ is the $(b,v,k,r)$-design $\Delta'$ obtained by interchanging
the roles of treatments and blocks, so the point-block incidence matrix
of $\Delta'$ is the transpose of that of $\Delta$. 
As the canonical efficiency factors of $\Delta'$ differ from 
those of $\Delta$ only in the number of times $1$ occurs, it follows that 
$\Delta$ is Schur optimal in a class $\cal C$ of $(v,b,r,k)$-designs 
if and only if $\Delta'$ is Schur optimal in the class of the duals
of the elements of $\cal C$  (see \citet{BC09}). 
We take the treatments of the \textit{dual} $S'$ of an $(n\times n)/k$ 
semi-Latin square $S$ to be the Cartesian product
$\{1,\ldots,n\}\times \{1,\ldots,n\}$, with treatment $(i,j)$ corresponding
to the $(i,j)$-cell of $S$, and then, for each treatment $\alpha$ of $S$,
the corresponding block in $S'$ is the set of those $(i,j)$ 
such that $\alpha$ is in the $(i,j)$-cell of $S$. In particular,
$S'$ is an $(n^2,nk,k,n)$-design. See \citet{Bailey11} 
for more on duals of semi-Latin squares, including applications.
Also relevant are parts of \citet{Suen82} and \citet{SC}. 

\citet{BMM95} studied and constructed affine resolvable designs and proved: 
\begin{thm}
\label{TH:AFFOPT}
Let $\Delta$ be an affine resolvable $(v,b,r,k)$-design with 
$r>2$, and let $s=v/k>1$. Then:
\begin{enumerate}
\item the canonical efficiency factors of $\Delta$
are $1-1/r$, with multiplicity $r(s-1)$, and $1$, with multiplicity $v-1-r(s-1)$;
\item
the affine resolvable $(v,b,r,k)$-designs
are precisely the Schur optimal designs in the class
of all resolvable $(v,b,r,k)$-designs.
\end{enumerate}
\end{thm}

We prove the following analogous result for uniform semi-Latin squares: 
\begin{thm}
\label{TH:USLSOPT}
Let $n>2$ and let $S$ be a uniform $(n\times n)/k$ semi-Latin square.
Then:
\begin{enumerate}
\item the canonical efficiency factors of $S$ are 
$1-1/(n-1)$, with multiplicity $(n-1)^2$, and $1$, with 
multiplicity $nk-1-(n-1)^2$;
\item
the uniform $(n\times n)/k$ semi-Latin squares
are precisely the Schur optimal designs in the class
of all $(n\times n)/k$ semi-Latin squares.
\end{enumerate}
\end{thm}

\begin{proof}
\citet{Soicher12} determined the canonical efficiency factors of $S$
and its Schur optimality. 
Here, we complete the proof of the theorem.

Let $T$ be any Schur optimal $(n\times n)/k$
semi-Latin square and let $T'$ be the dual of $T$. 
We shall show that the concurrence matrix $A'$ of $S'$
is equal to the concurrence matrix $B'$ of $T'$, 
showing that $T$ is uniform. 

The canonical efficiency factors of $S'$ are  
$1-1/(n-1)$, with multiplicity $(n-1)^2$ and $1$, with 
multiplicity $2n-2$. 
Now the argument in the proof of Theorem~3.4 of \citet{Soicher12}
shows that the 0-eigenspace of $A'$ 
is contained in the 0-eigenspace of $B'$,
so in particular $T'$ has at least $2n-2$ canonical efficiency factors
equal to 1. Then, since the sum of the canonical efficiency
factors is the same for $S'$ and $T'$, the Schur optimality of
both $S'$ and $T'$ implies that $T'$ must have precisely
$2n-2$ canonical efficiency factors equal to $1$, and the remaining ones 
equal to $1-1/(n-1)$. 

We now know that $A'$ and $B'$ have the same $nk$-eigenspace 
(spanned by the all-1 vector), as well as the same $0$-eigenspace.
It follows that the orthogonal complement of the direct sum of these 
two eigenspaces must be the eigenspace for the remaining
eigenvalue of both $A'$ and $B'$.
Thus $A'\mathbf{e}=B'\mathbf{e}$ as $\mathbf{e}$ 
runs over a basis of $\mathbb{R}^{n^2}$ 
(consisting of common eigenvectors), and so $A'=B'$.

Now concurrence in $S'$ and $T'$ is block intersection
size in $S$ and $T$, respectively, and so $S$ uniform implies $T$ uniform.
\end{proof}

\citet{Cal71} emphasised the importance of block designs with just two
distinct canonical efficiency factors, one of which is~$1$.  As a result, some
authors call this the C-property, and call such designs C-designs:  for
example, see \citet{Saha} and \citet{CKM}.

\section{Pairwise-variance aberration}

Now, given a collection ${\cal S}$ of Schur optimal designs in a given
class ${\cal C}$ of $(v,b,r,k)$-designs, we want a criterion to choose
between them.  For affine resolvable designs, \citet{Morgan10} proposed
choosing a design with minimum pairwise-variance (PV) aberration, a
combinatorial criterion which translates to a statistical one when the
Schur optimal $(v,b,r,k)$-designs under consideration are all connected
and all have the same two distinct canonical efficiency factors (and
no others).

\begin{defn}
Let $\Delta$ be a $(v,b,r,k)$-design. 
Define \[ \eta(\Delta):=(\eta_0(\Delta),\ldots,\eta_r(\Delta)),\] where
$\eta_i(\Delta)$ is the number of unordered pairs of distinct
treatments of $\Delta$ with concurrence equal to $i$. 
If $\Delta$ is connected and has at most two distinct 
canonical efficiency factors, then $\eta(\Delta)$ is
called the \textit{pairwise-variance (or PV) aberration} of $\Delta$.
Where $\Delta$ and $\Gamma$ are connected $(v,b,r,k)$-designs
having the same two distinct canonical efficiency factors,
the design $\Delta$ is considered to have \textit{smaller PV aberration} than
$\Gamma$ if $\eta(\Delta)$ is lexicograhically less than $\eta(\Gamma)$.
\end{defn}

Now the underlying $(nk,n^2,n,k)$-design of a uniform $(n\times n)/k$ semi-Latin square
is connected (since $n>2$) and has exactly two distinct canonical efficiency factors, 
$1-1/(n-1)$ and $1$. We thus make the following:

\begin{defn}
For $S$ a semi-Latin square with underlying block design $\Delta$,
we define $\eta(S)$ and $\eta_i(S)$ to be 
$\eta(\Delta)$ and $\eta_i(\Delta)$, respectively, and
if $S$ is uniform then we call $\eta(S)$ the   
\textit{PV aberration} of $S$. 
Where $S$ and $T$ are uniform $(n\times n)/k$ 
semi-Latin squares then $S$ is considered to have \textit{smaller PV aberration} than
$T$ if $\eta(S)$ is lexicograhically less than $\eta(T)$.
\end{defn}

These definitions are justified by the following result, which is a corollary
of Theorem~1 of \citet{Bailey09}.

\begin{thm}
Let $\Delta$ be a connected $(v,b,r,k)$-design having at most two distinct
canonical efficiency factors.  Then the variance of the estimator of
the difference of the effects of distinct treatments $\alpha$ and $\beta$ is a
function depending only on $r,k$, the canonical efficiency factors, and
the concurrence of $\alpha$ and $\beta$. 
Moreover, as a function of the concurrence of $\alpha$ and $\beta$
this function is strictly decreasing.
\end{thm}

Thus, by minimising PV aberration in an appropriate
class of designs, we are not only  
minimising the maximum pairwise-variance, but, 
for those designs in the class with the same maximum pairwise-variance, we are 
minimising the number of pairs of distinct treatments having that
maximum pairwise-variance, and when these numbers are the same, we are 
minimising the number of pairs with the next largest pairwise-variance, and so on.

\citet{Morgan10} studied affine resolvable designs
with minimum PV aberration, and placed an extensive catalogue
of these online at
\linebreak
\url{http://designtheory.org/database/v-r-k-ARD-MV/}.

We remark that when a $(v,b,r,k,\lambda)$-BIBD exists, with $b>v$,
it may be of interest to determine
the duals of $(v,b,r,k,\lambda)$-BIBDs with
minimum PV aberration, and for this, 
the block intersection size distribution of 
$(v,b,r,k,\lambda)$-BIBDs would need to be studied. 

We next present a general result providing designs with minimum PV aberration 
in certain circumstances, but first we need a definition.
For $s$ a positive integer, an \textit{$s$-fold inflation}
of a $(v,b,r,k)$-design or an $(n\times n)/k$ semi-Latin square is obtained by
replacing each treatment $\alpha$ by $s$
treatments $\sigma_{\alpha,1},\ldots,\sigma_{\alpha,s}$, such that
$\sigma_{\alpha,i}=\sigma_{\beta,j}$ if and only if $\alpha=\beta$
and $i=j$. In particular, an $s$-fold inflation of a $(v,b,r,k)$-design is
an $(sv,b,r,sk)$-design and an $s$-fold inflation of an $(n\times n)/k$ semi-Latin square
is an $(n\times n)/(sk)$ semi-Latin square. 

\begin{thm}
Suppose that $\Delta$ is a $(v,b,r,k)$-design, with $r>1$, such that
every pair of distinct non-disjoint blocks meet
in a positive constant number $\mu$ of treatments. Then
\begin{equation}
\eta_0(\Delta)\ge v(v-k-(r-1)(k-\mu))/2,
\label{eta0ge}
\end{equation}
with equality holding if and only if
\begin{equation}
\eta(\Delta)=(v(v-k-(r-1)(k-\mu))/2,vr(k-\mu)/2,0,\ldots,0,v(\mu-1)/2),
\label{EQ:eta}
\end{equation}
which happens if and only if $\Delta$ is a $\mu$-fold inflation of a
$(v/\mu,b,r,k/\mu)$-design with the property that every pair
of distinct non-disjoint blocks meet in just one treatment.
\label{TH:eta}
\end{thm}

\begin{proof}
Let $\alpha$ be any treatment of $\Delta$, and let $B_1,\ldots,B_r$
be the blocks containing $\alpha$ (in some fixed, but arbitrary,
order). We have that $|B_i\cap B_j|=\mu$ for $1\le i<j\le r$.

Now let $d_i$ be the number of treatments in $B_i$ that are not in
any of $B_1,\ldots,B_{i-1}$. Then $\alpha$ is concurrent with
exactly $d:=\sum_{i=1}^r d_i$ distinct
treatments (including $\alpha$ itself).
Now $d_1=k$, $d_2=k-\mu$, and for $i=3,\ldots,r,$
$d_i\le k-\mu$, with equality if and only if
$B_1,\ldots,B_i$ meet pairwise in exactly
the same $\mu$ treatments.
Thus $\alpha$ has concurrence $0$ with $v-d\ge v-k-(r-1)(k-\mu)$
treatments, with equality if and only if $B_1,\ldots,B_r$ meet pairwise in exactly
the same $\mu$ treatments.

From the preceding argument, we have that
the inequality ($\ref{eta0ge}$) holds, 
and that if equality holds then every treatment must have concurrence
$0$ with $v-k-(r-1)(k-\mu)$ treatments, concurrence
$1$ with $r(k-\mu)$ treatments, and 
concurrence $r$ with the remaining $\mu-1$ treatments other than itself,
in which case ($\ref{EQ:eta}$) holds.  

Now suppose that ($\ref{EQ:eta}$) holds,
and define a relation $\sim$ on the set of treatments by
$\alpha\sim \beta$ if and only if $\alpha$ and $\beta$
have concurrence $r$. Then $\sim$ is easily seen to
be an equivalence relation, with each equivalence class having
exactly $\mu$ elements. Choose equivalence class representatives
$\alpha_1,\ldots,\alpha_{v/\mu}$, and form the design
$\Gamma$ having these representatives as treatments
and blocks $B\cap \{\alpha_1,\ldots,\alpha_{v/\mu}\}$,
where $B$ runs over the blocks of $\Delta$. Then $\Gamma$
is a $(v/\mu,b,r,k/\mu)$-design such that every pair of distinct
non-disjoint blocks meet in just one treatment, and
$\Delta$ is a $\mu$-fold inflation of $\Gamma$.

Now if $\Delta$ is a $\mu$-fold inflation of 
a $(v/\mu,b,r,k/\mu)$-design with the property that 
that every pair of distinct non-disjoint blocks meet in just one treatment, 
then for each treatment $\alpha$ of $\Delta$, the blocks of $\Delta$
containing $\alpha$ meet in the same $\mu$ treatments,
and it follows that equality holds in ($\ref{eta0ge}$).
\end{proof}

We shall apply the preceding theorem to uniform semi-Latin squares.  
The \textit{superposition} of an $(n\times n)/k$ semi-Latin square
with an $(n\times n)/\ell$ semi-Latin square (with disjoint sets
of treatments) is made by superimposing the first square upon the
second, resulting in an $(n\times n)/(k+\ell)$ semi-Latin square.
For example, the uniform $(3\times 3)/4$ semi-Latin square 
($\ref{USLSexample}$) is a 2-fold inflation of a superposition of two MOLS of order~$3$. 
Indeed, it is easy to see that when $n>2$ and there exist $n-1$ 
MOLS of order $n$ that the 
superposition $T$ of these $n-1$ MOLS is a uniform
$(n\times n)/(n-1)$ semi-Latin square, and so, for every positive
integer $\mu$, a $\mu$-fold inflation of $T$ is a uniform $(n\times
n)/(\mu(n-1))$ semi-Latin square (see \citet{Soicher12}).

\begin{thm}
Suppose that $S$ is a uniform $(n\times n)/k$ semi-Latin square, 
with $\mu:=\mu(S)$. Then
\begin{equation*}
\eta_0(S)\ge nk^2/2, 
\end{equation*}
and the following are equivalent:
\begin{enumerate}
\item
$\eta_0(S)=nk^2/2$;
\item 
$\eta(S)=(nk^2/2,n^2k(k-\mu)/2,0,\ldots,0,nk(\mu-1)/2)$;
\item
$S$ is a $\mu$-fold inflation of a
superposition of $n-1$ MOLS of order $n$;
\item
$n-1$ MOLS of order $n$ exist and $S$ has minimum
PV aberration in the class of 
uniform $(n\times n)/k$ semi-Latin squares.
\end{enumerate}
\label{TH:USLSeta}
\end{thm}

\begin{proof}
The underlying block design of $S$ is an $(nk,n^2,n,k)$-design in which every pair
of distinct non-disjoint blocks meet in exactly $\mu=k/(n-1)$ treatments. Thus, by Theorem~\ref{TH:eta},
we have that $\eta_0(S)\ge nk(nk-k-(n-1)(k-k/(n-1)))/2=nk^2/2$.

Now we prove the equivalence of statements 1 to 4.

($1\Rightarrow 2$) This follows from Theorem~\ref{TH:eta}.

($2\Rightarrow 3$) Assume that statement~2 holds. It follows from Theorem~\ref{TH:eta} that
$S$ must be a $\mu$-fold inflation of a uniform $(n\times n)/(n-1)$ semi-Latin square $T$,
and, by Theorem~3.3 of \citet{Soicher12}, $T$ is a superposition of $n-1$ MOLS of order~$n$.

($3\Rightarrow 4$) Now assume that $S$ is a $\mu$-fold inflation of a
superposition of $n-1$ MOLS of order $n$. In particular, $n-1$ MOLS of order $n$ exist.
Now $\eta_0(S)=nk^2/2$, so if $U$ is any uniform $(n\times n)/k$ semi-Latin square 
then $\eta_0(U)\ge \eta_0(S)$, with equality 
if and only if $\eta(U)=\eta(S)$. Thus $S$ has minimum
PV aberration in the class of uniform $(n\times n)/k$ semi-Latin squares.

($4\Rightarrow 1$) Finally, assume that $n-1$ MOLS of order $n$ exist, and that 
$N$ is a $\mu$-fold inflation of the superposition of these $n-1$ MOLS.
Then $N$ is a uniform $(n\times n)/k$ semi-Latin square.
Now if $S$ has minimum
PV aberration in the class of uniform $(n\times n)/k$ semi-Latin squares, 
we must have $\eta_0(S)\le \eta_0(N)=nk^2/2$, but since 
$\eta_0(S)\ge nk^2/2$, we must have equality. 
\end{proof}

\begin{cor}
Let $n>2$ and assume that there exist $n-1$ MOLS of order $n$. 
Then for every positive integer $\mu$, the uniform $(n\times n)/(\mu(n-1))$
semi-Latin squares with minimum PV aberration are precisely the $\mu$-fold
inflations of the superpositions of $n-1$ MOLS of order $n$. 
\end{cor}

When the integer $n>2$ is a prime-power, a well-known construction of
\citet{Bose38} gives $n-1$ MOLS of order $n$, and so a $\mu$-fold
inflation of the superposition of these $n-1$ MOLS yields a uniform
(and hence Schur optimal) $(n\times n)/(\mu(n-1))$ semi-Latin square
with minimum PV aberration in the class of all Schur optimal (and
hence by Theorem~\ref{TH:USLSOPT} uniform) $(n\times n)/(\mu(n-1))$
semi-Latin squares.

For example, the $(5\times 5)/12$ uniform semi-Latin squares
were classified by \citet{Soicher13a}. Up to isomorphism, there are exactly 
277 such semi-Latin squares, and we calculated their PV aberrations. 
The least such PV aberration is \[(360, 1350, 0, 0, 0, 60),\] coming from a
3-fold inflation of the superposition of four MOLS of order~5. The next best 
PV aberration is \[(488, 1062, 128, 64, 0, 28),\]  and the worst is
\[(720, 450, 600, 0, 0, 0).\] This shows that when searching for designs
with minimum, or near minimum, PV aberration, one cannot restrict
the search to designs having concurrences differing by as little as
possible. This is contrary to the usual thinking for eigenvalue-based
optimality criteria, as discussed in \citet{JM77} and in Section 2.5 of
\citet{JW95}.

\section{The uniform $(6\times 6)/10$ semi-Latin squares} 
\label{SEC:classifying}

Let $n>2$. If $n$ is a prime power, we can use $n-1$ MOLS of order $n$ to construct 
a uniform $(n\times n)/(\mu(n-1))$ semi-Latin square with minimum PV aberration,
for all $\mu\ge 1$. This focuses attention on the case $n=6$.
There does not exist a uniform $(6\times 6)/5$ semi-Latin square,
since, by Theorem~3.3 of \citet{Soicher12}, such a square would
be a superposition of five MOLS of order~6, and even just two MOLS of order~6
do not exist. On the other hand, in Section~5 of \citet{Soicher12}, uniform $(6\times 6)/(5\mu)$
semi-Latin squares are constructed for all $\mu\ge 2$. 
 
In this section, we describe a complete classification of the uniform
$(6\times 6)/10$ semi-Latin squares, and display the 
one amongst these having least PV aberration. 
The computations described took place on a desktop PC with 16GB RAM and
an Intel(R) i7-6700 CPU running at 3.4GHz.

Consider now the Hamming graph $H(2,n)$. This graph has vertex-set
$\{1,\ldots,n\}\times \{1,\ldots,n\}$, with distinct vertices $(i,j)$
and $(i',j')$ joined by an edge if and only if $i=i'$ or $j=j'$. 
We observe that a block design $\Delta$ is the dual of a uniform 
$(n\times n)/(\mu(n-1))$ semi-Latin square if and only if: 
\begin{itemize}
\item
the treatments of $\Delta$ are the vertices of $H(2,n)$; 
\item
each block of $\Delta$ is a co-clique (independent set)
of size $n$ of $H(2,n)$ (and has multiplicity 
at most $\mu$ in $\Delta$);
\item
the concurrence of distinct treatments $\alpha,\beta$ of $\Delta$ 
is $0$ or $\mu$, according as $\{\alpha,\beta\}$ is an edge or non-edge of $H(2,n)$.
\end{itemize}
Note that, in particular, these conditions imply that every treatment of 
$\Delta$ is in exactly $\mu(n-1)$ blocks. 

The (graph) automorphism group $\mathrm{Aut}(H(2,n))$ of $H(2,n)$
is the wreath product $S_n\wr C_2$, which is generated by the direct
product $D:=S_n\times S_n$ of two copies of the symmetric group on
$\{1,\ldots,n\}$ and an element $\tau$ satisfying $\tau^2=1$ and
$\tau(g,h)=(h,g)\tau$ for all $(g,h)\in D$. If $(i,j)$ is a vertex of
$H(2,n)$ then the image of $(i,j)$ under $(g,h)\in D$ is $(i^g,j^h)$
and the image of $(i,j)$ under $\tau$ is $(j,i)$.

Now suppose that $S$ is a uniform $(n\times n)/(\mu(n-1))$ semi-Latin
square, and $\xi\in \mathrm{Aut}(H(2,n))$. Then the image of the dual $S'$ of
$S$ under $\xi$ is the block design whose treatments are the $\xi$-images of
the treatments of $S'$ and whose blocks are obtained by applying $\xi$ to 
every treatment in every block of $S'$. This image is the dual of
a uniform $(n\times n)/(\mu(n-1))$ semi-Latin square. 
Now let $T$ be any uniform $(n\times n)/(\mu(n-1))$ semi-Latin square.
By Theorem~3 of \citet{Soicher13b}, we have that $S$ and $T$
are isomorphic as semi-Latin squares if and only if there is an element
of $\mathrm{Aut}(H(2,n))$ mapping the dual $S'$ of $S$ to the dual $T'$
of $T$. Moreover, due to the structure of both $S'$ and $T'$ as duals
of uniform semi-Latin squares, any block design isomorphism from $S'$
to $T'$ must be a graph automorphism of $H(2,n)$, and so we have that $S$
and $T$ are isomorphic as semi-Latin squares if and only if $S'$ and $T'$
are isomorphic as block designs.

We classify the uniform $(6\times 6)/10$ semi-Latin squares via
backtrack searches (see, for example, Section~6.2 of \citet{GO07},
for an introduction to using backtrack search for the enumeration of
block designs).  Our searches are for the block multisets of the duals
of the uniform semi-Latin squares we seek to classify.  We represent a
block multiset as a set of (block,multiplicity)-pairs, where the blocks in
these pairs are distinct and the associated multiplicity of a block gives
the number of times that block occurs in the block multiset.  Our backtrack
searches work (block,multiplicity)-pair by (block,multiplicity)-pair.

We define a \textit{partial solution} to be a set of
(block,multiplicity)-pairs, such that the blocks are distinct co-cliques
of size $6$ of $H(2,6)$, the (block) multiplicities 
are 1 or 2, and no non-edge of
$H(2,6)$ is contained in more than two blocks (counting multiplicities). A
\textit{solution} is a partial solution for which every non-edge of
$H(2,6)$ is contained in exactly two blocks (counting multiplicities).
Hence, the solutions are precisely the block multisets of the 
duals of the uniform $(6\times 6)/10$ semi-Latin squares. 

We first programmed a partial backtrack search exploiting the 
automorphism group of $H(2,6)$, using the \textsf{GAP} system \citep{GAP}
and adapting code from its \textsf{DESIGN} \citep{DESIGN} and \textsf{GRAPE}
\citep{GRAPE} packages. This search was used to generate a sequence
\[ (P_1,A_1),(P_2,A_2),\ldots,(P_t,A_t),\] 
where each $P_i$ is a partial solution and its corresponding $A_i$
is a set of (block,multiplicity)-pairs, such that no block of $A_i$ is
a block of $P_i$, and the following hold:
\begin{itemize}
\item
each isomorphism class of duals of
uniform $(6\times 6)/10$ semi-Latin squares has at least one representative
whose block multiset is a solution consisting of some $P_i$ extended 
by elements belonging to the corresponding $A_i$; 
\item
each $P_i$ (considered as a multiset of independent sets of $H(2,6)$) 
has trivial stabiliser in $\mathrm{Aut}(H(2,6))$.
\end{itemize}
For our program, it turned out that $t=2214$, and the search 
ran in under four minutes. 

Then, for each pair $(P_i,A_i)$, we used a newly developed
C~program to perform a backtrack search to determine all the solutions
which are extensions of $P_i$ by elements from $A_i$. 
The total run time for this step was about eight and 
a half hours, or on average, about 14 seconds for each $i$. 

Finally, we took all the (1340930 as it turned out) solutions found by the
C~program backtrack searches and determined isomorphism class representatives
amongst all the duals of uniform $(6\times 6)/10$ semi-Latin squares
having those solutions as their block multisets. We did this using
the \textsf{DESIGN} package making heavy use of the \textsf{bliss} program
\citep{bliss} via \textsf{GRAPE}. The run time for this step was
a little over five hours. 

We found that, up to isomorphism, there are exactly 8615 uniform $(6\times
6)/10$ semi-Latin squares. These semi-Latin squares, their PV aberrations,
and their duals are available from
\url{http://www.maths.qmul.ac.uk/~lsoicher/usls/}.
There is a unique (up to isomorphism) uniform $(6\times 6)/10$ semi-Latin
square $M$ with least PV aberration, which is 
\[ (532,906,294,30,6,0,2).\]
We present this $M$ in Figure~\ref{FIG:1}. The automorphism group of
the dual of $M$ has order 12, but the automorphism group of $M$ has
order 48, since there are automorphisms of $M$ fixing every cell, but
interchanging the treatments in one or both of the pairs of treatments
with concurrence~6. We note that $\eta_0(M)=532$, which is well off the 
lower bound of $300$ given by Theorem~\ref{TH:USLSeta}. 

\begin{figure}
\begin{footnotesize}
\[
\addtolength{\arraycolsep}{-0.6\arraycolsep}
\begin{array}{|rrrrr|rrrrr|rrrrr|rrrrr|rrrrr|rrrrr|}
  \hline
  1&2 &3 & 4 & 5 & 11 & 12 & 13 & 14 & 15 & 21 & 22 & 23 & 24 & 25
  & 31 & 32 & 33 & 34 & 35 & 41 & 42 & 43 & 44 & 45 & 51 & 52 & 53 & 54 & 55\\
  6 & 7 & 8 &9 & 10 & 16 & 17 & 18 & 19 & 20 & 26 & 27 & 28 & 29 & 30
  & 36 & 37 & 38 & 39 & 40 & 46 & 47 & 48 & 49 & 50 & 56 & 57 & 58 & 59 & 60\\
  \hline
  11 & 12 & 21 & 22 & 31 & 1 & 2 & 23 & 24 & 33 & 3 & 4 & 13 & 14 & 35 &
  5 & 6 & 15 & 16 & 25 & 7 & 8 & 17 & 18 & 27 & 9 & 10 & 19 & 20 & 29\\ 
  32 & 41 & 42 & 51 & 52 & 34 & 43 & 44 & 53 & 54 & 36 & 45 & 46 & 55 & 56
  & 26 & 47 & 48 & 57 & 58 & 28 & 37 & 38 & 59 & 60 & 30 & 39 & 40 & 49 & 50\\
  \hline
  17 & 19 & 25 & 27 & 35 & 7 & 9 & 21 & 22 & 36 & 1 & 2 & 15 & 20 & 31 &
  3 & 10 & 13 & 18 & 23 & 4 & 5 & 11 & 16 & 29 & 6 & 8 & 12 & 14 & 24\\
  39 & 43 & 45 & 53 & 57 & 40 & 47 & 48 & 55 & 59 & 37 & 41 & 49 & 58 & 60
  & 28 & 44 & 50 & 51 & 52 & 30 & 32 & 33 & 54 & 56 & 26 & 34 & 38 & 42 & 46\\
  \hline
  14 & 15 & 24 & 28 & 37 & 3 & 5 &26 & 27 & 32 & 8 & 9 & 11 & 19 & 33 &
  1 & 2 & 12 & 17 & 29 & 6 & 10 & 13 & 20 & 21 & 4 & 7 & 16 & 18 & 23\\
  40 & 47 & 50 & 54 & 56 & 39 & 46 & 49 & 51 & 60 & 38 & 44 & 48 & 52 & 57
  & 30 & 42 & 45 & 55 & 59 & 22 & 34 & 35 & 53 & 58 & 25 & 31 & 36 & 41 & 43\\
  \hline
  13 & 16 & 23 & 29 & 34 & 4 & 8 & 28 & 30 & 31 & 6 & 10 & 12 & 18 & 32 &
  7 & 9 & 11 & 20 & 24 & 1 & 2 & 14 & 19 & 25 & 3 & 5 & 15 & 17 & 21\\
  38 & 48 & 49 & 55 & 60 & 35 & 42 & 50 & 57 & 58 & 39 & 43 & 47 & 54 & 59 &
  27 & 41 & 46 & 53 & 56 & 26 & 36 & 40 & 51 & 52 & 22 & 33 & 37 & 44 & 45\\
  \hline
  18 & 20 & 26 & 30 & 33 & 6 & 10 & 25 & 29 & 37 & 5 & 7 & 16 & 17 & 34 &
  4 & 8 & 14 & 19 & 21 & 3 & 9 & 12 & 15 & 23 & 1 & 2 & 11 & 13 & 27\\
  36 & 44 & 46 & 58 & 59 & 38 & 41 & 45 & 52 & 56 & 40 & 42 & 50 & 51 & 53 &
  22 & 43 & 49 & 54 & 60 & 24 & 31 & 39 & 55 & 57 & 28 & 32 & 35 & 47 & 48\\
  \hline
\end{array}
\]
\end{footnotesize}
\caption{Uniform $(6\times 6)/10$ semi-Latin square with least PV aberration}
\label{FIG:1}
\end{figure}

The next best PV aberration of a uniform $(6\times 6)/10$ 
semi-Latin square is 
\[(532,912,276,48,0,0,2),\] and the worst is \[(600,720,450,0,0,0,0).\] 

Using the \textsf{DESIGN} package, we found that no dual of a uniform
$(6\times 6)/10$ semi-Latin square is resolvable. Equivalently, no uniform
$(6\times 6)/10$ semi-Latin square is a superposition of Latin squares,
which answers a question raised by \citet{Soicher13a}.

We have checked that the results of our computations are consistent
with some previous partial classifications of uniform $(6\times 6)/10$
semi-Latin squares done by the authors using the \textsf{DESIGN} package.
These include finding that (up to
isomorphism) there are exactly 5828 uniform $(6\times 6)/10$ semi-Latin
squares whose dual has a non-trivial automorphism, exactly 7 uniform
$(6\times 6)/10$ semi-Latin squares having at least two treatments with
concurrence~6, and (see \citet{Soicher13a}) exactly 98 uniform $(6\times
6)/10$ semi-Latin squares having no concurrence greater than 2.
We have also checked that our programs agree with \citet{Soicher13a}
that, up to isomorphism, there are exactly 10 uniform $(5\times 5)/8$
semi-Latin squares and exactly 277 uniform $(5\times 5)/12$ semi-Latin 
squares. 

\section{A new construction of uniform semi-Latin squares}

Suppose that $n>2$ and there exist $n-1$ mutually orthogonal Latin squares $\Lambda_1$,
\ldots, $\Lambda_{n-1}$ of order~$n$ with disjoint sets of symbols $\mathcal{L}_1$, \ldots,
$\mathcal{L}_{n-1}$.  We present a construction of an
$((n+1)\times (n+1))/(\mu n)$ uniform semi-Latin square with $\mu=n-2$, 
which generalises the uniform $(6\times 6)/15$ semi-Latin square used in the 
proof of Theorem~5.1 of \citet{Soicher12}.

For $i=1$, \ldots, $n-2$ and $j=1$,\ldots, $n$, 
put $\mathcal{L}_{ij} = \{(\alpha,j): \alpha\in\mathcal{L}_i\}$.
Put $\bar{\mathcal{L}}_i = \mathcal{L}_{i1} \cup \cdots \cup \mathcal{L}_{in}$.
Create an $n$-fold inflation $\bar{\Lambda}_i$ of $\Lambda_i$ using the symbols in
$\bar{\mathcal{L}}_i$.

Put $\bar{\mathcal{L}}_{n-1} = \mathcal{L}_{n-1} \times \{1, \ldots, n-2 \}$, and create an
$(n-2)$-fold inflation $\bar{\Lambda}_{n-1}$ of $\Lambda_{n-1}$ using the symbols in
$\bar{\mathcal{L}}_{n-1}$.

Superpose $\bar{\Lambda}_1$, \ldots, $\bar{\Lambda}_{n-1}$ to give an
$(n \times n)/ \ell$ 
semi-Latin square $S$, where $\ell = (n-2)n + n-2 = (n-2)(n+1)$.

Add an extra row and an extra column to this.  For $i=1$, \ldots, $n$, $j=1$, \ldots, $n$ and
$t=1$, \ldots, $n-2$, each cell $(i,j)$  contains a unique symbol from ${\mathcal{L}}_{tj}$:
remove this from cell $(i,j)$ and insert it in cells $(i,n+1)$ and $(n+1,j)$.  Put all the symbols in
$\bar{\mathcal{L}}_{n-1}$ into cell $(n+1,n+1)$.

Now we have an array $\bar{S}$ of size $((n+1) \times (n+1))/k$ where
$k=\ell - (n-2) = n(n-2) = \left| \bar{\mathcal{L}}_{n-1}\right|$.
Every symbol in 
$\bar{\mathcal{L}}_1 \cup \cdots \cup \bar{\mathcal{L}}_{n-2} \cup
\bar{\mathcal{L}}_{n-1}$ 
occurs precisely once in each row and once in each column, so
$\bar{S}$ is a semi-Latin square. 

For $1\leq i,j \leq n$, cell $(i,j)$ contains symbols $(\alpha,1)$, \ldots, $(\alpha,n-2)$
for a single symbol $\alpha$ in $\mathcal{L}_{n-1}$: thus it has $n-2$ symbols in common
with cell $(n+1,n+1)$.

Suppose that $i\ne i'$ and $j\ne j'$ with $i$, $i'$, $j$, $j'$ in $\{1, \ldots,n\}$.  
Then there is exactly one value
of $t$ such that the original cells $(i,j)$ and $(i',j')$ have the same symbol in $\Lambda_t$.
If $t=n-1$ then the symbols in common to the new cells $(i,j)$ and $(i',j')$ in
$\bar{S}$ are precisely those in $\bar{\mathcal{L}}_{n-1}$: there are
$n-2$ of these symbols.  If $1\leq t\leq n-2$ then the symbols in common to the
new cells $(i,j)$ and $(i',j')$ are precisely those in 
$(\bar{\mathcal{L}}_{t} \setminus \mathcal{L}_{tj}) \setminus \mathcal{L}_{tj'}$: there are
$n-2$ of these.

We have shown that new cell $(i,j)$ in $\bar{S}$
has $n-2$ symbols in common with new cell $(i',j')$ for all $j'$ in
$\{1, \ldots, n\} \setminus \{j\}$.  This accounts for $(n-1)(n-2)$ symbols in new
cell $(i,j)$.  Since it has no symbols in common with new cell $(i',j)$, the remaining 
$n-2$ symbols in new cell $(i,j)$ must be in cell $(i',n+1)$.
Similarly, new cell $(i,j)$ has $n-2$ symbols in common with cell $(n+1,j')$.

For $1\leq i\leq n$ and $1 \leq j\leq n$, we have shown that  cell $(i,n+1)$ has precisely 
$n-2$ symbols in common with new cell $(i',j)$ when $i'\ne i$ and none in common with new
cell $(i,j)$.  For any fixed $j$, this accounts for $(n-1)(n-2)$ symbols in column $j$.
Hence the remaining $n-2$ symbols in cell $(i,n+1)$ must be in cell $(n+1,j)$.

This completes the proof that $\bar{S}$ is a uniform semi-Latin square.

\begin{thm}
\label{TH:etabarS}
  Let $\bar{S}$ be the $((n+1) \times(n+1))/((n-2)n)$ uniform
  semi-Latin square constructed above.  If $n\geq 5$ then $n-2>2$ and
  \begin{eqnarray*}
    \eta_{n+1}(\bar{S}) & = & \frac{n(n-2)(n-3)}{2},\\
    \eta_{n-2}(\bar{S}) & = & \frac{n^2(n-1)(n-2)}{2},\\
    \eta_{2}(\bar{S}) & = & \frac{n^2(n-2)(n-3)(2n-1)}{2},\\
    \eta_{1}(\bar{S}) & = & \frac{n(n-1)(n-2)(n^3-4n^2 +8n -2)}{2},\\
    \eta_{0}(\bar{S}) & = & \frac{n^2(n-2)(3n^2-9n+4)}{2},
  \end{eqnarray*}
and $\eta_m(\bar{S}) =0$ for all other non-negative integers~$m$.
\end{thm}

\begin{proof}
  Suppose that $1\leq i \leq n-2$ and $\alpha \in \mathcal{L}_i$.
  For $j$ and $j'$ in $\{1, \ldots, n\}$ with $j \ne j'$, the concurrence of
  $(\alpha,j)$ and $(\alpha,j')$ is $n-2$.
  If $\beta \in \mathcal{L}_i \setminus
  \{\alpha\}$, then the concurrence of $(\alpha,j)$ and $(\beta,j)$ is $1$
  and there is one other value of $j'$ such that the concurrence of
  $(\alpha,j)$ and $(\beta,j')$ is $1$; otherwise this concurrence is~$0$.

  Suppose that $1 \leq i' \leq n-2$ and $i' \ne i$.  Then we can show that,
  of the $n^2$ elements in $\bar{\mathcal{L}}_{i'}$, $2n-1$ have
  concurrence~$2$ with $(\alpha,j)$, $(n-1)(n-3)$ have concurrence~$1$ with
  $(\alpha,j)$, and the remaining $2(n-1)$ have concurrence~$0$ with
  $(\alpha,j)$.

  Now consider $\gamma$ in $\mathcal{L}_{n-1}$.  For $j$ and $j'$ in
  $\{1, \ldots, n-2\}$ with $j\ne j'$,
  the concurrence of $(\gamma,j)$ and $(\gamma,j')$ is  $n+1$.
  If $\delta \in \mathcal{L}_{n-1} \setminus \{\gamma\}$ and
  $j'\in \{1, \ldots, n-2\}$ then the concurrence of $(\gamma,j)$ and
  $(\delta,j')$ is~$1$. If $\alpha\in\mathcal{L}_i$, as above, then the
  concurrence of $(\gamma,j)$ and $(\alpha,j')$ is $0$ for one value
  of $j'$ and is $1$ for the remaining $n-1$ values of $j'$.

  It follows that $\eta_{n+1}(\bar{S}) = n(n-2)(n-3)/2$, and
  that  $\eta_{n-2}(\bar{S}) = n(n-2)n(n-1)/2$.  Concurrence~$2$
  occurs only between $\bar{\mathcal{L}}_i$ and $\bar{\mathcal{L}}_{i'}$ for
  $i\ne i'$, so we find that $\eta_2(\bar{S}) = (n-2)(n-3)/2 \times
  n^2(2n-1)$.

  Concurrence~$1$ occurs $(n-2)n^2 \times 2(n-1)/2$ times within treatment
  sets $\bar{\mathcal{L}}_1$, \ldots, $\bar{\mathcal{L}}_{n-2}$, and
  $(n-2)(n-3)n^2/2 \times (n-1)(n-3)$ times between such sets.  It occurs
  $n(n-1)/2 \times (n-2)^2$ times within $\bar{\mathcal{L}}_{n-1}$,
  and $n(n-2) \times (n-2) \times n(n-1)$ times between $\bar{\mathcal{L}}_{n-1}$
  and other treatments.
  Finally, concurrence~$0$ occurs $(n-2)n^2 \times (n-1)(n-2)/2$ times within
  $\bar{\mathcal{L}}_1$, \ldots, $\bar{\mathcal{L}}_{n-2}$,
  $(n-2)(n-3)n^2/2 \times 2(n-1)$ times between these sets, and
  $n(n-2) \times (n-2) \times n$ times between these and $\mathcal{L}_{n-1}$.
\end{proof}

Let $n\ge 5$. We observe that the lower bound for $\eta_0$ 
given by Theorem~\ref{TH:USLSeta} applied to $\bar{S}$ is
$(n+1)n^2(n-2)^2/2$, which is
$(n+1)(n-2)/(3n^2-9n+4)$ times $\eta_0(\bar{S})$.
Thus, when $n+1$ is a prime power, $\bar{S}$ is far from
optimal with respect to PV aberration. 
However, when $n+1$ is not a prime power, we do not know how
far from optimal $\bar{S}$ is with respect to 
PV aberration. 

Putting $n=5$ in Theorem~\ref{TH:etabarS} 
shows that the uniform $(6 \times 6)/15$ semi-Latin square
constructed by \citet{Soicher12} has PV aberration
\[
(1275, 1890, 675, 150, 0, 0, 15).
\]
On the other hand, the least PV aberration of any uniform $(6\times 6)/15$
semi-Latin square found by the authors so far is
\[
( 1260, 1943, 630, 125, 40, 0, 7 ),
\]
and a uniform $(6\times 6)/15$ semi-Latin square with this
PV aberration is available from
\url{http://www.maths.qmul.ac.uk/~lsoicher/usls/}. 

Suppose now $n>2$ is a prime power. Then there exist $n-1$ MOLS of order
$n$, and so we can make a uniform $((n+1)\times (n+1))/((n-2)n)$
semi-Latin square $\bar{S}$ as above.  
In this case, there is also the construction
of Theorem~4.3 of \citet{Soicher12} 
(equivalent to the construction of Corollary~4.1.2 of \citet{Suen82}) 
which gives a uniform semi-Latin
square $T$ of size $((n+1) \times (n+1))/((n-1)n)$ (and also one of size
$((n+1) \times (n+1))/((n-1)n/2)$ when $n$ is odd).  Now let $s$ and $t$
be non-negative integers, not both zero. Let $U$ be an $s$-fold inflation
of $\bar{S}$ if $t=0$, let $U$ be a $t$-fold inflation of $T$ if
$s=0$, and otherwise let $U$ be the superposition of an $s$-fold
inflation of $\bar{S}$ and a $t$-fold inflation of $T$. Then $U$
is a uniform $((n+1)\times (n+1))/(\mu n)$ semi-Latin square with
$\mu=s(n-2)+t(n-1)$.  Since $n-2$ and $n-1$ are coprime, 
by the well-known solution of the ``Frobenius coin problem" for
two denominations, every integer greater than or equal to $(n-3)(n-2)$
is a non-negative integer linear combination of $n-2$ and $n-1$.
Therefore, when $n>2$ is a prime power, there exists a uniform semi-Latin
square of size $((n+1) \times (n+1))/(\mu n)$ for every positive
integer $\mu\geq (n-3)(n-2)$ (and for every positive integer $\mu\ge
(n-3)((n-1)/2-1)$ if in addition $n$ is odd).

\section{Constructing affine resolvable designs and BIBDs from uniform semi-Latin squares}
\label{SEC:constructing}

Let $n$ be any integer greater than~2, 
let $S$ be a uniform $(n \times n)/(\mu(n-1))$ semi-Latin square, with rows
$R_1,\ldots,R_n$, and let
$\Delta(S)$ be its underlying block design. Obtain the block design 
$\Delta_1(S)$ from $\Delta(S)$ by adding, for each $i=1,\ldots,n$, 
$\mu$ new treatments $R_{i,1},\ldots,R_{i,\mu}$, 
each incident precisely with the blocks in row $R_i$. In
$\Delta_1(S)$, the set of blocks in each column form a replicate, and every pair
of blocks in different columns have exactly $\mu$ treatments in common.  Hence 
$\Delta_1(S)$ is an affine resolvable $(\mu n^2,n^2,n,\mu n)$-design.  The
analogous construction using columns in place of rows gives another affine resolvable design
$\Delta_2(S)$, which may or may not be isomorphic to $\Delta_1(S)$. 
However, if $T$ is any uniform $(n\times n)/(\mu(n-1))$ semi-Latin
square isomorphic to $S$, then the isomorphism classes of $\Delta_1(T)$
and $\Delta_2(T)$ are those of $\Delta_1(S)$ and $\Delta_2(S)$,
in some order.

Each of the $8615$ uniform $(6\times 6)/10$ semi-Latin squares $S$
classified in Section~\ref{SEC:classifying} gives two affine resolvable
designs $\Delta_1(S)$ and $\Delta_2(S)$. We found that, up to block
design isomorphism (determined by the \textsf{DESIGN} package using
\textsf{bliss} via \textsf{GRAPE}), these give in total $16875$ affine
resolvable $(72,36,6,12)$-designs. The only isomorphisms of affine
resolvable designs resulted from the $355$ uniform $(6\times 6)/10$
semi-Latin squares $T$ having an automorphism interchanging rows and
columns, in which case $\Delta_1(T)$ is isomorphic to $\Delta_2(T)$.

An affine resolvable $(72,36,6,12)$-design is equivalent to an orthogonal
array $\mathrm{OA}(72,6,6,2)$.  Recent work finding all isomorphism
classes of orthogonal arrays for certain parameter tuples, such as that
by \citet{BuluMar}, \citet{BuluRy1,BuluRy2}, \citet{GBR} and \citet{SEN},
does not yet include the orthogonal arrays $\mathrm{OA}(72,6,6,2)$, but
our present work can be used to produce $16875$ pairwise non-isomorphic
orthogonal arrays having this parameter tuple.
An arbitrary orthogonal array $\mathrm{OA}(72,6,6,2)$ $A$ is isomorphic to
one of these $16875$ orthogonal arrays if and only if $A$ is isomorphic to
an orthogonal array $B$ with the property that for each $i=1,2,3,4,5,6$,
there are two rows of $B$ of the form $(i,i,i,i,i,i)$. An instance of
this, which extends to an $\mathrm{OA}(72,7,6,2)$,
is given in Example~4.1.2 of \citet{Suen82}.

For every uniform $(n\times n)/(\mu(n-1))$ semi-Latin square $S$ 
we can make a third design $\Delta_3(S)$ from $\Delta(S)$ by adding
$\mu$ new treatments for each row and $\mu$ new treatments for each column, 
as above (so each new treatment is incident with the blocks in its
corresponding row or column), and then taking the dual. 
This gives a block design with 
$n^2$~treatments in $\mu n(n-1) + 2\mu n$ blocks of
size~$n$, in which every pair of distinct treatments has concurrence~$\mu$.  In other words,
$\Delta_3(S)$ is an $(n^2, \mu n(n+1), \mu(n+1),n, \mu)$-BIBD.

When $\mu=1$, this construction is essentially that of \citet{Bose38} 
to obtain an affine plane of order $n$ (i.e{.} an 
$(n^2,n(n+1),n+1,n,1)$-BIBD) from $n-1$ MOLS of order $n$. 
When $\mu\geq 2$, the BIBDs we construct have repeated blocks.  
However, their parameter tuples seem to have no overlap with those of
the BIBDs with repeated blocks constructed by \citet{RABPJC}.

If $S$ and $T$ are isomorphic uniform $(n\times n)/k$ semi-Latin squares
then $\Delta_3(S)$ and $\Delta_3(T)$ must be isomorphic, but the converse
need not hold. For example, it was shown by \citet{OP} (and confirmed
by \citet{EW}) that there are exactly 15 sets of eight MOLS of order~9,
up to trisotopism, that is, up to uniform permutation of the rows, uniform
permutation of the columns, permuting the symbols separately within each
Latin square, and optionally, simultaneously transposing the Latin 
squares. It follows, from Theorem~3.3 of \citet{Soicher12}, 
that there are exactly 15 isomorphism classes of uniform $(9\times 9)/8$
semi-Latin squares. However, there are only seven isomorphism classes
of affine planes of order~9 (see, for example, \citet{OP}), and so there
are non-isomorphic uniform $(9\times 9)/8$ semi-Latin squares $S$ and $T$
such that $\Delta_3(S)$ and $\Delta_3(T)$ are isomorphic affine planes.

We have, however, verified (via the \textsf{DESIGN} package using
\textsf{bliss} via \textsf{GRAPE}), that the $8615$ $(36,84,14,6,2)$-BIBDs
$\Delta_3(S)$ obtained from the $8615$ uniform $(6 \times 6)/10$
semi-Latin squares $S$ classified in Section~\ref{SEC:classifying} are
pairwise non-isomorphic.  According to \citet{MR07}, only five BIBDs
with these parameters were known at that time.

\section*{Acknowledgement} 
The support from EPSRC grant EP/M022641/1 
(CoDiMa: a Collaborative
Computational Project in the area of Computational Discrete Mathematics)
is gratefully acknowledged.

\end{document}